\documentclass[12pt]{amsart}
\usepackage{amsmath}

\def\CC{{\rm\kern.24em\vrule width.02em height1.4ex depth-.05ex\kern-.26em C}}
\def\ZZ{{\rm\kern.26em\vrule width.02em height0.5ex depth0ex\kern.04em\vrule width.02em height1.47ex depth-1ex\kern-.34em Z}}
\def\BB{{\rm\kern.24em\vrule width.02em height1.4ex depth-.05ex\kern-.26em B}}

\newcommand {\BOX} {\rule{2mm}{2mm} \bigskip}
\newcommand{\ra}{\rightarrow}
\newcommand {\eps} {\epsilon}
\newcommand {\ol} {\overline}
\newcommand {\dist}{\text{\rm dist}}

\newtheorem{lem}{Lemma}[section]
\newtheorem{theo}[lem]{Theorem}

\newtheorem{remark}[lem]{Remark}

\numberwithin{equation}{section}

\begin{document}
\title{An embedding of $\CC$ in $\CC^2$ with hyperbolic complement}
\author{Gregery T. Buzzard and John Erik Forn\ae ss}
\date{}

%\tableofcontents

\begin{abstract}
Let $X$ be a closed, $1$-dimensional, complex subvariety of $\CC^2$
and let $\ol{\BB}$ be a closed ball in $\CC^2 - X$.  Then there exists a
Fatou-Bieberbach domain $\Omega$ with $X \subseteq \Omega \subseteq
\CC^2 - \ol{\BB}$ and a biholomorphic map $\Phi: \Omega \ra \CC^2$
such that $\CC^2 - \Phi(X)$ is Kobayashi hyperbolic.  As corollaries,
there is an embedding of the plane in $\CC^2$ whose complement is
hyperbolic, and there is a nontrivial Fatou-Bieberbach domain containing any
finite collection of complex lines.
\end{abstract}

\renewcommand{\thefootnote}{}
\footnote{Research at MSRI is 
   supported in part by NSF grant DMS-9022140}

\maketitle

%%%%%%%%%%%%%%%%%%%%%%%%%%%%%%%%%%%%%%%%%%%%%%%%%%%%%%%%%%%%%%%%%%%%%%
\section{Introduction}
This paper was motivated by the question of whether there exists a
proper holomorphic embedding $\Phi
:\CC \ra \CC^2$ such that $\CC^2 - \Phi(\CC)$ is Kobayashi hyperbolic.

For the case of a polynomial embedding $P:\CC \ra \CC^2$, Abhyankar
and Moh \cite{am} and Suzuki \cite{suzuki}, proved that there exists a
polynomial automorphism $F$ of $\CC^2$ such that $(F \circ
P)(\CC) = \CC \times \{0\}$.  In contrast, Forstneric, Globevnik and
Rosay \cite{fgr} have shown that there exists a proper holomorphic
embedding $\Phi: \CC \ra \CC^2$ such that for no automorphism $H$ of
$\CC^2$ is it true that $(H \circ \Phi)(\CC) = \CC \times \{0\}$.
This leaves open the question of whether the complement of an embedded
copy of $\CC$ is always Kobayashi hyperbolic.  We answer this
question in the following theorem.

\begin{theo}   \label{thm:main}
There exists a proper holomorphic embedding $\Phi:\CC \ra \CC^2$ such
that $\CC^2 - \Phi(\CC)$ is Kobayashi hyperbolic.
\end{theo}

This theorem is a corollary of theorem~\ref{thm:universal}.  As
additional corollaries, we obtain the following theorems.

\begin{theo}
There exists a proper holomorphic embedding of $\CC$ into $\CC^2$ such
that any nonconstant holomorphic image of $\CC$ intersects this
embedding infinitely many times.
\end{theo}

\begin{theo}
Let $X$ be any finite set of complex lines in $\CC^2$.  Then there
exists a Fatou-Bieberbach domain $\Omega$ containing $X$ with
$\ol{\Omega} \neq \CC^2$. 
\end{theo}

This last theorem answers a question of Rosay and Rudin \cite{rr}, who
constructed a Fatou-Bieberbach domain containing the coordinate
axes and asked if a Fatou-Bieberbach domain can contain more than 2
complex lines. 

%%%%%%%%%%%%%%%%%%%%%%%%%%%%%%%%%%%%%%%%%%%%%%%%%%%%%%%%%%%%%%%%%%%%%%
\section{Controlling images of the unit disk}

We first prove a lemma showing that certain collections of linear
disks give control on certain (nonlinear) maps from the unit disk to
$\CC^2$. 

Given two concentric balls in $\CC^2$, we find a collection of linear
disks near the boundary of the inner ball such that any map from the
unit disk into the larger ball which maps $0$ near the center of the
two balls must map most of the disk into the smaller ball.  We make
this more precise in the following lemma.  The proof
uses a normal families argument and is similar in spirit to the
construction of non-tame sets in \cite{rr}.

For notation, let $\pi_j$ denote projection to the $j$th coordinate,
$j=1,2$, let $\Delta(0,r)$ be the (open) disk of radius $r$ centered at $0$
in $\CC$, and let $\BB(0,r)$ denote the (open) ball of radius $r$ centered at
the origin in $\CC^2$.

\begin{lem}  \label{lemma:disks}
Let $k \in \ZZ^+$, $0 \leq n_1<n_2< n_3$. Let $X$ 
be a closed, complex, $1$-dimensional subvariety of $\CC^2$.
Then there exist 
finitely many affine complex linear maps $L_j: \CC \ra \CC^2$ with
$\ol{\Delta_j} := L_j(\ol{\Delta(0,1)}) \subseteq \BB(0,n_2)-
(\ol{\BB(0,n_1)} \cup X)$ having pairwise disjoint closures
such that if 
\begin{align*}
\phi:\Delta(0,1) & \ra \BB(0,n_3) - \cup_j \Delta_j,
\\
\phi(0) & \in \ol{\BB(0,n_1/2)},
\end{align*}
with $\dist(\phi(0),X) \geq 1/k$,
then $\phi(\Delta(0,1-1/2^k)) \subseteq \BB(0, n_2)$.
\end{lem}

In fact we prove a stronger result allowing small perturbations in
place of the linear disks obtained above.

\begin{lem}   \label{lemma:perturb}
With the hypotheses of lemma~\ref{lemma:disks}, there exist disks
$\Delta_j = L_j(\Delta(0,1))$ as in that lemma, plus $\delta>0$,
such that if $\Delta_j'$ is the image of a
holomorphic map which is within $\delta$ of $L_j$ on $\Delta(0,1)$ for
all $j$, then for any
\begin{align*}
\phi:\Delta(0,1) & \ra \BB(0,n_3) - \cup_j \Delta_j',
\\
\phi(0) & \in \ol{\BB(0,n_1/2)},
\end{align*}
with $\dist(\phi(0),X) \geq 1/k$,
we have $\phi(\Delta(0,1-1/2^k)) \subseteq \BB(0, n_2)$.
\end{lem}

{\bf Proof of lemma~\ref{lemma:perturb}:}  Let $n_1<r<n_2$, and let
$\{a_j\}$ and $\{b_j\}$ be countable dense subsets 
of $\BB(0,r) - \ol{\BB(0,n_1)}$ and $\BB(0,n_2) -
\ol{\BB(0,r)}$, respectively, such that $\pi_l a_j \neq \pi_l a_m$
if $j \neq m$, $l=1,2$, and likewise for $b_j$.

For each $j$, let $A_j(z) =
(\alpha_j z,0) + a_j$ and $B_j(w) = (0,\beta_j w) + b_j$, where
$\alpha_j>0$ is chosen maximal such that $A_j(\Delta(0,2)) \subseteq
\BB(0,r) - (\ol{\BB(0,n_1)} \cup X)$, and analogously for
$\beta_j$ with disks in the outer shell.

For notation, we let
$A_{j,m}'$ and $B_{j,m}'$ denote functions holomorphic and within $1/m$ of
$A_j$ and $B_j$ on $\Delta(0,1)$.  For such functions, let
$$D_m := \cup_{j=1}^m (A_{j,m}'(\Delta(0,1)) \cup
B_{j,m}'(\Delta(0,1))).$$

To reach a contradiction, assume that for each $m$ there exists $D_m$
as above and 
$\phi_m:\Delta(0,1)  \ra \BB(0,n_3) - D_m $, $\phi_m(0) \in
\ol{\BB(0,n_1/2)}$ with $\dist(\phi_m(0), X) \geq 1/k$, and
$\phi(\Delta(0,1-1/2^k)) \not \subseteq \BB(0, n_2)$. 
 
The set $\{\phi_m\}$ is a normal family since the image of each
$\phi_m$ lies inside a fixed ball, so we may assume that some 
subsequence converges to a map $\phi:\Delta(0,1)  \ra \BB(0,n_3)$
with $\phi(0) \in \ol{\BB(0,n_1/2)}$, 
$\dist(\phi(0), X) \geq 1/k$, and $\phi(\Delta(0,1)) \not
\subseteq \BB(0,n_2)$.  
Then there exist nonempty open sets $\Omega_1, \Omega_2 \subseteq
\Delta(0,1)$ such that $\phi(\Omega_1) \subseteq \BB(0,r) -
\ol{\BB(0, n_1)}$ and $\phi(\Omega_2) \subseteq \BB(0,n_2) -
\ol{\BB(0, r)}$.  

We claim that $\pi_2 \phi$ is constant.  If not, then since $\phi(0)
\notin X$, there exists $z_1
\in \Omega_1$ such that $\phi(z_1) \notin X$ and $(\pi_2 \phi)'(z_1)
\neq 0$, and there exists a subsequence $a_{j_l} \ra \phi(z_1)$.  Moreover,
there exists $c>0$ such that $\alpha_{j_l} \geq c$ for all $l$.  

By continuity, there exists a neighborhood $V$ of $z_1$ such that
$|\pi_1 \phi(z) - \pi_1\phi(z_1)|<c/2$ for all $z \in V$ and such that
$\pi_2 \phi(V)$ is a neighborhood of $\pi_2 \phi(z_1)$.  For $l$
large, $\pi_2 a_{j_l} \in \pi_2 \phi(V)$, and hence by choice of $V$
there exists $z_1' \in V$ such that $\phi(z_1') \in
A_{j_l}(\Delta(0,1))$.  But then $\phi(z) - A_{j_l}(\zeta)$ has an
isolated zero in $\Delta^2(0,1)$, which persists under small
perturbations.  Hence $\phi_m$ must intersect $D_m$ for large $m$, a
contradiction, so $\pi_2 \phi \equiv const$.  A similar argument
for $\pi_1 \phi$ shows that $\phi \equiv const$, which contradicts
$\phi(0) \in \ol{\BB(0,n_1/2)}$ and $\phi(\Delta(0,1)) \not \subseteq
\BB(0,n_2)$.   

Thus, no such sequence $\{\phi_m\}$ exists, so we obtain the lemma by
taking $\{L_j\}$ to be some finite subset of $\{A_j\} \cup \{B_j\}$
and $\delta$ sufficiently small.  $\BOX$  \bigskip

%%%%%%%%%%%%%%%%%%%%%%%%%%%%%%%%%%%%%%%%%%%%%%%%%%%%%%%%%%%%%%%%%%%%%%
\section{Approximation of linear disks by complex subvarieties}

In the following lemma, we use the notation $reg(X)$ to denote
the set of regular points of a subvariety.

\begin{lem}  \label{lemma:variety}
Let $K \subseteq \CC^2$ be polynomially convex and let $X$ be a
closed, $1$-dimensional, complex subvariety of $\CC^2$.  Let $L
\subseteq X$ be compact such that for each $p \in \text{reg}(X)-(K 
\cup L)$ there exists a curve $\eta:[0,1) \ra
\text{reg}(X)-(K \cup L)$ with $\eta(0)=p$ and $\lim_{t \ra
1} \|\eta(t) \| = \infty$.  

Then $K \cup L$ is polynomially convex.
\end{lem}

{\bf Proof:}  Although the ideas in the proof are standard, we include
it for convenience.  

Let $p \in \CC^2 - (K \cup L)$.  If $p \not \in X$, then there
exists an entire $g$ such that $g \equiv 0$ on $X$ but $g(p) \neq 0$.
Also, there exists $f$ entire such that $|f|<1$ on $K$ but
$|f(p)|>1$.  Then for some $m$, we have $|f^m g(p)|>1$ and $|f^m g|<1$
on $K \cup L$, and we can approximate $f^m g$ by a polynomial with the same
properties. 

If $p \in X$, then the preceding argument shows that if $K_p := K \cup
L \cup \{p\}$, then $\hat{K}_p \subseteq K \cup X$, where $\hat{K}_p$
is the polynomial hull of $K_p$.  The hypotheses on $L$ imply that $X
\cap K_p$ is Runge in $X$, so we can approximate functions holomorphic
on $X \cap K_p$ by functions holomorphic on $X$.

Let $f(p)=2$ and $f\equiv 1/2$ on $(X \cap K_p)- \{p\}$, and
approximate $f$ by a function $g$ holomorphic on $X$ such that
$|f-g|<1/4$ on the 
domain of $f$.  By \cite[theorem 18, ch. VIII]{gunning-rossi}, we may
extend $g$ to be holomorphic on $\CC^2$, and then we may restrict $g$
to a  neighborhood $V$ of $X$ such that 
$|g|<1$ on $V \cap K$.  Composing with a convex function
$\sigma$, we obtain a psh function $\sigma \circ |g|$ which is $0$ on
a neighborhood of $(X \cap K_p)-\{p\}$ and larger than $1$ at $p$.  We can
extend this function by $0$ so that it is psh in a neighborhood of $K \cup
X$, $0$ on $K \cup L$, and larger than $1$ at $p$, then restrict it to
a Runge neighborhood of $\hat{K}_p$.  

Using standard $\ol{\partial}$ techniques, it follows that there exists a
polynomial $P$ on $\CC^2$ such that $|P(p)|>1$ and $|P|\leq 1$ on $K
\cup L$, so $K \cup L$ is polynomially convex.  $\BOX$  \bigskip

In the following lemma, we start with a closed, $1$-dimensional, complex
subvariety $X$ of $\CC^2$, an automorphism $\Psi$, and some linear
disks contained outside 
$\Psi(X)$ and outside the ball of radius $n$.  We construct an automorphism
of $\CC^2$ such that the image of $\Psi(X)$ under this automorphism contains
pieces approximating each of the linear disks, and such that this
automorphism is near the identity on $\BB(0,n)$ and on some large
piece of $\Psi(X)$.  We do this in such a way that
given some polynomially convex set contained in a ball disjoint from
$\Psi(X)$ and $\BB(0,n+1)$, the image of this set is disjoint from
$\BB(0,n+2)$.  

\bigskip
For notation, let $\eps>0$, $R, n \geq 0$, and let $\Psi$ be an
automorphism of $\CC^2$.  Let $X$ be a closed, $1$-dimensional, complex
subvariety of $\CC^2$, and let $L_j: \CC \ra \CC^2$, $j=1, \ldots, N$
be finitely many affine complex linear maps with 
$\Delta_j := L_j(\Delta(0,1))$ having pairwise disjoint closures such that
$\ol{\Delta_j} \cap \Psi(X) = \emptyset$ and $\ol{\Delta_j} \cap
\ol{\BB(0,n)} = \emptyset$ for all $j$.  Let $\BB_0= \BB(p_0,r_0)$ be
a ball with $\ol{\BB_0} \cap (\ol{\BB(0,n+1)} \cup \Psi(X)\cup \cup_j
\ol{\Delta_j}) = \emptyset$, and let $K_0 \subseteq \BB_0$ be
polynomially convex and compact.

\begin{lem}   \label{lemma:approx}
There exists $H:\CC^2 \ra \CC^2$ an automorphism 
such that $\|H-I\|< \eps$ on $\BB(0,n) \cup \Psi(\BB(0,R) \cap X)$
and such that for each $j$, there is a
submanifold $\sigma_j$ of $H \Psi(X)$ such that the orthogonal projection
of $\sigma_j$ onto $\Delta_j$ is a diffeomorphism and such that
$\sigma_j$ is near $\Delta_j$ in the sense of
lemma~\ref{lemma:perturb}.  Moreover, there exists a ball $\BB_0'$
with $\ol{\BB_0'} \cap (\ol{\BB(0,n+2)} \cup H\Psi(X)) = \emptyset$
and $H(K_0) \subseteq \BB_0'$. 
\end{lem}

{\bf Proof:}  
Choose $p_1, \ldots, p_N \in \Psi(reg(X) - \ol{\BB(0,R)})-
\ol{\BB(0,n)}$  and choose pairwise disjoint $C^2$ curves
$\gamma_j:[0,1] \ra \CC^2$, $j = 0, \ldots, N$ which
are disjoint from $\Psi(X \cap \ol{\BB(0,R)}) \cup \ol{\BB(0,n)}$ such
that $\gamma_j(0) = p_j$ and such that if $j \geq 1$, then
$\gamma_j(1)= L_j(0)$, while $\|\gamma_0(1)\| > n+2$.  Let $W_j$ be a
neighborhood of $\gamma_j$ such that the $W_j$ have pairwise disjoint
closures which are also disjoint from $\Psi(X \cap \ol{\BB(0,R)}) \cup
\ol{\BB(0,n)}$, and such that $\ol{\Delta_j} \subseteq W_j$ if
$j \geq 1$ and $\ol{\BB(p_0,r_0)} \subseteq W_0$.
For $j \geq 1$, let $r_j>0$ such that $\Psi(X) \cap \ol{\BB(p_j,r_j)}
\subseteq W_j$ and such that there is an embedding $g_j: \Delta(0,2)
\ra \CC^2$ with $g_j(\ol{\Delta(0,1)}) = \Psi(X) \cap \ol{\BB(p_j,r_j)}$. 

Let $K =  \Psi(X \cap \ol{\BB(0,R)}) \cup \ol{\BB(0,n)} \cup
\ol{\BB_0} \cup (\cup_{j=1}^N g_j(\ol{\Delta(0,1)}))$.  In order to
apply the approximation result in \cite[theorem 2.1]{fr},
we will construct a family $\Phi_t$ of maps which are biholomorphic
in a neighborhood of $K$ and $C^2$ in $t$ with $\Phi_0 = I$.  We will
construct $\Phi_t$ so that
$\Phi_1$ is the identity on $\BB(0,n) \cup \Psi(X \cap \BB(0,R))$ and 
maps each $g_j(\Delta(0,1))$ to a disk near $\Delta_j$, and
such that $\Phi_t(K)$ is polynomially convex for all $t$, then
approximate $\Phi_1$ by a global automorphism.

We show how to construct $\Phi_t$ around a neighborhood of $p_1$.  For
simplicity, assume $p_1=0$ and $g_1(0)=p_1$.  Now, given any complex
linear map $L$ such that $L(\ol{\Delta(0,1)}) \subseteq W_1$ is
tangent to $g_1(\Delta(0,1))$ at $p_1$,
we can first use the family $g_1(\chi_1(t) g_1^{-1}(p))$
with $\chi_1(0)=1$, $\chi_1(1)$ small and positive, to
shrink $g_1(\Delta(0,1))$ to a nearly linear disk tangent to
$L(\Delta(0,1))$, then use the family $\Lambda_1(t) g_1(\chi_1(1)
g_1^{-1}(p))$ to expand the small disk to approximate
$L(\Delta(0,1))$.  

Note that by choosing $\chi_1(1)$ very small, we can
make the diameter of $g_1(\chi_1(1) \Delta(0,1))$ as small as we like.
In particular, we can translate this disk along $\gamma_1$ and use a
one-parameter family of rotations to make the image tangent to
$\Delta_j$, then expand as before to make the new image approximate
$\Delta_j$.  Note that all of this can be done within $W_j$, and that
at each stage, the image of $\Delta(0,1)$ is contained in an
affine linear image of $X$.  Moreover, we can make the family
$C^2$ in $t$.  Reparametrizing, and using a similar construction for
each $p_j$, we can define $\Phi_t$ on $\cup_{j=1}^N
g_1(\ol{\Delta(z_j,r_j)})$.  We can use a similar procedure to shrink
$\ol{\BB_0}$ and move it along $\gamma_0$ inside $W_0$ so that the
image is $\ol{\BB_0'}$ as in the statement.

Finally, each $g_j$, $j \geq 1$ can be extended to be a biholomorphic
embedding of $\ol{\Delta^2(0,1)}$ into $W_j$, so we can use the same
argument to define $\Phi_t$ on a neighborhood of $\cup_{j=1}^N
g_j(\ol{\Delta(0,1)})$, and similarly extend it to a neighborhood of
$\ol{\BB_0}$.  Define $\Phi_t \equiv I$ on a neighborhood of  
$\Psi(X \cap \ol{\BB(0,R)}) \cup \ol{\BB(0,n)}$ for all $t$.

Since the union of two disjoint closed balls is polynomially convex,
and since the remaining hypotheses of the previous lemma are satisfied
for all $t$, we see that $\Phi_t(K)$ is polynomially convex for each
$t$.  Hence by \cite[theorem 2.1]{fr}, there is a neighborhood $V$ of $K$
such that $\Phi_1$ can be approximated uniformly on $V$ by
automorphisms $H$ of $\CC^2$.

Choosing $\Phi_t$ such that $\Phi_1 g_j(\Delta(0,1))$ is close to
$\Delta_j$ for all $j$ and $\Phi_1(K_0) \subseteq \BB_0'$, then
choosing an automorphism $H$ close to $\Phi_1$, we obtain the lemma.
$\BOX$ \bigskip

%%%%%%%%%%%%%%%%%%%%%%%%%%%%%%%%%%%%%%%%%%%%%%%%%%%%%%%%%%%%%%%%%%%%%%
\section{Main theorem}

\begin{theo}  \label{thm:universal}
Let $X$ be a closed, $1$-dimensional, complex subvariety of $\CC^2$ and
$\BB_0$ a ball with $\ol{\BB_0} \cap X = \emptyset$.  Then there
exists a domain $\Omega \subseteq \CC^2 - \ol{\BB_0}$ containing $X$
and a biholomorphic map $\Phi$ from $\Omega$ onto $\CC^2$ such that
$\CC^2 - \Phi(X)$ is Kobayashi hyperbolic.  Moreover, all nonconstant
images of $\CC$ in $\CC^2$ intersect $\Phi(X)$ in infinitely many
points.
\end{theo}

\begin{remark}
Since $\Phi$ in the statement of the theorem is biholomorphic, it is
necessarily a proper holomorphic embedding of $X$ into $\CC^2$.
\end{remark}

{\bf Proof:}
The first step is to construct $\Phi$ as the limit of automorphisms of
$\CC^2$ which are constructed inductively.  We will define a sequence
$\Phi_m$ of automorphisms, numbers $R_m \nearrow \infty$, and finitely
many affine complex linear disks $\Delta_j^m$ such that

\begin{itemize}
\item[$(1_m)$] $\|\Phi_m(p)\| \geq m+1$ if $ p \in X$ and $\|p\| \geq R_m$,
\item[$(2_m)$] $\|\Phi_m(p)\| \geq m+1$ if $ p \in \ol{\BB_0}$, 
\item[$(3_m)$] $\|\Phi_{m+1} - \Phi_m\| \leq 1/2^m$ on $X \cap
	\ol{\BB(0,R_m)}$,
\item[$(4_m)$] $\|\Phi_{m+1} \circ \Phi_m^{-1} - I\| \leq 1/2^m$ on
	$\ol{\BB(0,m)}$ and
\item[$(5_m)$] if $l < m$ and $\phi:\Delta(0,1) \ra
\BB(0,l+2)-\cup_j(\Delta_j^l)'$ as in lemma~\ref{lemma:perturb} with
$\phi(0) \in \ol{\BB(0,l/2)}$ and 
dist$(\phi(0), \Phi_l(X)) \geq 1/l$, then $\phi(\Delta(0,1-1/2^l))
\subseteq \BB(0,l+1)$.
\end{itemize}

Changing coordinates by a translation, we may assume $\ol{\BB_0} \cap
\ol{\BB(0,1)} = \emptyset$, so we can 
take $\Phi_0 = I$, $R_0 = 1$, and choose a ball $\BB_0^0$ such that
$\ol{\BB_0} \subseteq \BB_0^0$ and $\ol{\BB_0^0} \cap (\ol{\BB(0,1)}
\cup X) = \emptyset$.  

For the inductive construction,
suppose we have $m \geq 0$, $\Phi_l$ and $R_l$ for $l \leq
m$ and linear discs $\Delta^l_j$ for $l<m$ satisfying $(1_l)$ and
$(2_l)$ for $l \leq m$ and $(3_l)$, $(4_l)$ and $(5_l)$ for $l<m$.
Suppose also that there exist balls $\BB^l_0$, $l\leq m$, with
$\Phi_l(\ol{\BB_0}) \subseteq \BB_0^l$ and $\ol{\BB_0^l} \cap
(\ol{\BB(0,l+1)} \cup \Phi_l(X)) = \emptyset$, $l \leq m.$  Note that 
$\Phi_m(\ol{\BB_0})$ is polynomially convex since $\Phi_m$ is a global
automorphism.  Suppose also that $\Phi_m(X \cap \BB(0,R_m))$ has
submanifolds approximating each $\Delta_j^l$ as in lemma~\ref{lemma:perturb}.

\bigskip

We use lemma~\ref{lemma:perturb} with $m$, $m+1$ and $m+2$ 
in place of $n_1$, $n_2$ and $n_3$, respectively, with $\Phi_m(X)$ in
place of $X$, and with $k=m$.  This gives 
finitely many affine complex linear $L_j^m: \CC \ra \CC^2$
with $\ol{\Delta_j^m}:= L_j^m(\ol{\Delta(0,1)}) \subset
\BB(0,m+1) - \ol{\BB(0, m)}$ having pairwise disjoint closures
such that if $\phi:\Delta(0,1) \ra \BB(0, m+2) - \cup_j
(\Delta_j^m)'$ as in lemma~\ref{lemma:perturb}, and $\phi(0)
\in \ol{\BB(0, m/2)}$ with dist$(\phi(0), \Phi_m(X)) \geq 1/m,$
then $\phi(\Delta(0,1-1/2^m))\subseteq \BB(0,m+1).$ 

\bigskip

We then apply lemma~\ref{lemma:approx} with $R=R_m$, $n=m$, $\Psi =
\Phi_m$, $X$ unchanged, $\{L_j^m\}$ in place of $\{L_j\}$, $\BB_0^m$
in place of $\BB_0$, and $K_0 = \Phi_m(\ol{\BB_0})$.  This gives an
automorphism $H: \CC^2 \ra \CC^2$ such that $\|H-I\| \leq \eps$ on 
$\ol{\BB(0,m)} \cup \Phi_m(X \cap \ol{\BB(0,R_m)})$ and such that there
are submanifolds $\sigma_j^m$ of $H\phi_m(X)$ which are
close enough to the $\Delta_j^m$ for lemma~\ref{lemma:perturb} and
such that there is a ball $\BB_0^{m+1}$ with $\ol{\BB_0^{m+1}} \cap
(\ol{\BB(0,m+2)} \cup H \Phi_m (X)) = \emptyset$ and $H
\Phi_m(\ol{\BB_0}) \subseteq \BB_0^{m+1}$.  

By taking $\eps$ sufficiently small, 
$(1_m)$ together with the fact that $\ol{\Delta_j^l} \subseteq
\BB(0,m)$ for $l<m$ implies that we
can make $H\Phi_m(X)$ pass close to all of these
disks, in the sense of lemma~\ref{lemma:perturb}.  Thus,  
taking $\Phi_{m+1} = H \Phi_m$ and $R_{m+1} \geq R_m +1$ large, we
obtain the inductive condition on submanifolds of $\Phi_{m+1}(X)$, as
well as $(1_{m+1}), (2_{m+1}), (3_m), (4_m)$ and $(5_m)$.

\bigskip
To show that $\{\Phi_m\}$ converges to give a map $\Phi:\Omega \ra \CC^2$,
let $\Omega_m:= \Phi_m^{-1}(\BB(0,m))$.  If $p \in \Omega_m,$
then $\| \Phi_{m+1}(p)-\Phi_m(p) \|= \|(\Phi_{m+1}\circ
\Phi_m^{-1})(\Phi_m(p))-I(\Phi_m(p))\|< \frac{1}{2^m}.$  Hence
$\|\Phi_{m+1}(p)\| < m+1/2^m< m+1,$ so $p\in \Omega_{m+1}.$  It follows
that the sequence $\{\Phi_m\}$ converges locally uniformly on $\Omega:=\cup
\Omega_m$ to a map $\Phi.$  \bigskip 

Suppose $p \in \Omega_{m+1} - \Omega_m.$ Then 

\begin{align*}
\|\Phi(p)- \Phi_{m+1}(p)\| & =  \| \lim_{k \ra \infty} (\Phi_k(p)-\Phi_{m+1}(p))\| \\  
                           & \leq \sum_{k \geq m+1}\|\Phi_{k+1}(p)-\Phi_k(p)\|\\ 
                           & < \frac{1}{2^m}
\end{align*}

Hence $\Phi(p) \in \BB(0,m+1+1/2^m) - \BB(0,m-1/2^m).$  It follows
that $\Phi$ is a proper map from $\Omega$ to $\CC^2.$ In particular,
$\Phi$ must have maximal rank $2$. Since $\Phi$ is a limit of
automorphisms, it follows that $\Phi$ has nonvanishing Jacobian
everywhere in $\Omega.$  A standard argument implies that $\Phi$ is
injective on $\Omega$, and $(4_m)$ implies that $\BB(0,m-1) \subseteq
\Phi_k(\Omega_m)$ for $k \geq m$, so that $\Phi$ is surjective.  Hence
$\Phi$ is a biholomorphism from $\Omega$ onto $\CC^2$.  \bigskip

{}From $(2_m)$ we see that $\ol{\BB_0} \cap \Omega_m= \emptyset$ for all
$m,$ hence $\ol{\BB_0}\cap \Omega= \emptyset.$  \bigskip

If $p \in X,$ then $p \in X\cap \ol{\BB(0,R_m)}$ for some $m$. It
follows from $(3_m)$ that $\|\Phi_k(p)\| \leq
\|\Phi_m(p)\|+\sum_{j=m}^{k-1}1/2^j < \|\Phi_m(p)\|+1$ for all $k>m.$
Hence if $k$ is large enough, $\|\Phi_k(p)\|<k$ and so $p \in
\Omega_k.$ Therefore $X \subset \Omega.$ \bigskip

Before continuing, note that the approximation of $\Delta_j^k$ by each
$\Phi_m(X \cap \BB(0,R_m))$ implies that $\Phi(X)$ has submanifolds which
approximate each $\Delta_j^k$ in the sense of lemma~\ref{lemma:perturb}.

\bigskip

We show next that $\CC^2 - \Phi(X)$ is
Kobayashi hyperbolic.  If not, there is a point $p$ in 
$\CC^2-\Phi(X)$, a nonzero tangent vector $\xi \in T_p \CC^2$,
and a sequence of holomorphic maps $\phi_k:\Delta(0,k+1) \ra \CC^2 -
\Phi(X)$ such that $\phi_k(0) = p_k$, $\phi_k'(0)= \xi_k$, and
$(p_k,\xi_k) \ra (p,\xi).$ 

Since $p \notin \Phi(X)$, there is an integer $m_0>1$ such that
$p \in \BB(0,m_0/2)$ and dist$(p,\Phi(X))> 1/m_0$.
From $(1_m)$ and $(3_m)$, there exist $k_1>1$, $m_1 > m_0$ such
that dist$(p_k, \Phi_m(X)) \geq 1/m_0$ for $k\geq k_1$, $m\geq m_1$.

Fixing $k\geq k_1$, we see that $\phi_k(\Delta(0,k)) \subseteq \BB(0,m+2)$
for some $m \geq m_1$ large.  Since $\phi_k$ misses $\Phi(X)$, we see
that $\psi_k(z):= \phi_k(z/k)$ maps
$\Delta(0,1)$ into $\BB(0,m+2)- \cup_j(\Delta_j^m)'$, $\psi_k(0)= p_k$
and dist$(\psi_k(0), \Phi_m(X)) \geq 1/m$.  Thus
$\psi_k(\Delta(0,1-1/2^m)) \subset \BB(0,m+1)$ by
$(5_m)$.  By induction, we
see that $\psi_k(\Delta(0,\Pi_{m\geq m_1}(1-1/2^m))) \subset
\BB(0,m_1+1)$.  Hence there exists $r>0$ independent of $k$ such
that $\|\psi_k(z)\| < m_1+1$ for all $|z|< r$.
This implies that $\|\phi_k(z)\| < m_1+1$
for all $|z|< kr$ and hence that $\phi_k'(0) \ra 0$, a contradiction.
Thus, $\CC^2-\Phi(X)$ is Kobayashi hyperbolic.

Finally, if $g: \CC \ra \CC^2$ is a nonconstant holomorphic map which
intersects $\Phi(X)$ only finitely many times, then these points of
intersection are contained in some large ball.  By reparametrizing, we
may assume that $g(0) \notin \Phi(X)$, then use an
argument like the one just given to show that the image of $g$ must be
contained in some large ball, a contradiction.
$\BOX$ \bigskip

The theorems in section 1 now follow immediately by taking $X$ to be the
$z$-axis for the first two theorems and to be any finite collection of
complex lines for the third theorem. In Theorem 1.3 one can replace $X$ by
any countable union of closed subvarieties of $\CC^2$ avoiding a fixed ball.
The proof of Theorem 4.1 still goes through. Moreover, if $X$ is a dense, countable
union of closed subvarieties of $\CC^2$ and $S$ is any discrete
set of points in the complement of $X$, then there exists a Fatou-Bieberbach
domain $\Omega, \; X \subset \Omega \subset \CC^2 \setminus S.$

Note that theorem~\ref{thm:universal} can be applied to any Riemann
surface which admits a proper holomorphic embedding into $\CC^2$, so
in particular, there is an embedding of the disk whose complement is
Kobayashi hyperbolic.

%%%%%%%%%%%%%%%%%%%%%%%%%%%%%%%%%%%%%%%%%%%%%%%%%%%%%%%%%%%%%%%%%%%%%%

\bigskip

\small
\noindent
Gregery T. Buzzard\\
Department of Mathematics\\
The University of Michigan\\
Ann Arbor, MI 48109, USA\\
and \\
MSRI\\
1000 Centennial Drive\\
Berkeley, CA 94720, USA\\

\noindent
John Erik Forn\ae ss\\
Department of Mathematics\\
The University of Michigan\\
Ann Arbor, MI 48109, USA\\

\end{document}